\documentclass{amsart}

\usepackage{hyperref}
\usepackage{fullpage, pstricks, subfigure, multirow}

\title{On Lengths of Rainbow Cycles}
\author{Boris Alexeev}
\email{borisa@mit.edu}
\thanks{This research is supported by MIT's Paul E. Gray Endowed Fund for UROP}
\date{November 8, 2005}
\subjclass[2000]{05C15}
\keywords{Complete graph, coloring, rainbow cycle, colorful cycle}

\theoremstyle{plain}  % default
\newtheorem {theorem}            {Theorem}
\newtheorem {lemma}     [theorem]{Lemma}
\newtheorem {claim}     [theorem]{Claim}
\newtheorem {corollary} [theorem]{Corollary}
\theoremstyle{definition}
\newtheorem*{definition*}        {Definition}
\newtheorem {problem}            {Problem}
\newtheorem {question}  [theorem]{Question}
\newtheorem*{conjecture*}        {Conjecture}
\theoremstyle{remark}
\newtheorem*{remark*}            {Remark}

\newcommand{\N}{\ensuremath{\mathbb{N}}}
\newcommand{\Z}{\ensuremath{\mathbb{Z}}}
\newcommand{\upto}{\ensuremath{\nearrow}}
\newcommand{\dnto}{\ensuremath{\searrow}}

\DeclareMathOperator{\col}{col}

\begin{document}
  \begin{abstract}
    We prove several results regarding edge-colored complete graphs and rainbow cycles, cycles with
    no color appearing on more than one edge.  We settle a question posed by Ball, Pultr, and
    Vojt\v{e}chovsk\'{y} \cite{ball:rainbow} by showing that if such a coloring does not contain a
    rainbow cycle of length $n$, where $n$ is odd, then it also does not contain a rainbow cycle of
    length $m$ for all $m$ greater than $2n^2$.  In addition, we present two examples which
    demonstrate that this result does not hold for even $n$.  Finally, we state several open
    problems in the area.
  \end{abstract}
\maketitle

\section{Introduction}

A rainbow cycle within an edge-colored graph is a cycle all of whose edges are colored with distinct
colors.  Rainbow cycles, sometimes called colorful or totally multicolored cycles in other sources,
have been introduced in many different contexts.  For example, Burr, Erd\H os, S\'os, Frankl, and
Graham \cite{MR1000076,MR1170778} studied which graphs have the property that regardless of how
their edges are $r$-colored (for a fixed integer $r$), the coloring contains a rainbow subgraph, in
particular a rainbow cycle of a certain length.

From another perspective, Erd\H os, Simonovits, and S\'os \cite{MR0379258} investigated a function
$f(n, C_k)$, defined as the maximum number of colors in which the edges of the complete graph $K_n$
on $n$ vertices may be colored so that the coloring contains no rainbow $k$-cycles; they conjectured
that $f(n, C_k) = n\cdot\bigl((k - 2)/2 + 1/(k - 1)\bigr) + O(1)$ for $n \ge k \ge 3$ and proved the
case $k=3$.  Alon \cite{MR693025} proved the conjecture for $k = 4$ and derived an upper bound for
general $k$; Jiang and West \cite{MR2037072} further improved these bounds and mentioned that the
conjecture has been proven for all $k \le 7$.  Finally, Montellano-Ballesteros and Neumann-Lara
\cite{MR2190794} recently proved the conjecture completely.  More research has occurred in related
areas; these references are not intended as a comprehensive survey of the area, but rather only a
small sample.

Within this paper, we build on the research of Ball, Pultr, and Vojt\v{e}chovsk\'{y}
\cite{ball:rainbow}, who studied rainbow cycles within edge-colored complete graphs.  In particular,
they asked when the existence of a rainbow cycle of a certain length forces the existence of a
rainbow cycle of another length, approaching the problem with applications for distributive
lattices.  To begin, it is easy to see that within such a coloring, the absence of rainbow
$n$-cycles and rainbow $m$-cycles implies the absence of rainbow $(n+m-2)$-cycles.  (For example, by
this fact alone, the absence of rainbow $3$-cycles implies there are no rainbow cycles at all, while
the absence of rainbow $4$-cycles implies there are no even-length rainbow cycles.  In general, if
there are no rainbow $n$-cycles, then there are no cycles with lengths congruent to $2\bmod (n-2)$.)
Ball et~al. asked whether or not the restrictions implied by this observation were the only
restrictions on the lengths of rainbow cycles in a coloring.

We answer this question in the negative by showing that there are more, and much stronger,
restrictions on the possible lengths.  In particular, our main result states that if a coloring does
not contain a rainbow cycle of length $n$, where $n$ is odd, then it also does not contain a rainbow
cycle of length $m$ for all $m \ge 2n^2$.  These results follow from the observation above combined
with intermediate results that show that the absence of a rainbow $n$-cycle (again, with $n$ odd)
implies the absence of rainbow cycles of lengths $\binom{n}{2}$ and $3n-6$.

However, the absence of a rainbow cycle of even length does not put the same restrictions on the
lengths of longer rainbow cycles.  For example, Ball et~al. showed that there are colorings with no
even-length rainbow cycles, but with rainbow cycles of all odd lengths up to the number of vertices
in the graph.  We present these colorings, as well as original colorings in which there are rainbow
cycles of all lengths except those congruent to $2\bmod 4$.

We briefly describe the organization of the paper.  We begin with preliminaries, including
definitions and the proof of the observation above, allowing us to formulate the original question
as well as the main result formally.  In order to provide a flavor of colorings that avoid rainbow
cycles, we present the ``even-length'' examples next.  We then prove that the absence of a rainbow
$n$-cycle (for $n$ odd) implies the absence of rainbow $\binom{n}{2}$-cycles, and show how this
prohibits all sufficiently large cycles (with lengths on the order of a cubic in $n$).  Next, we
show the stronger result that the absence of a rainbow $n$-cycle (for $n$ odd) also prohibits
rainbow $(3n-6)$-cycles, implying the absence of all sufficiently large cycles (with lengths on the
order of a quadratic in $n$, in particular $2n^2$, the main result).  At the end of the paper, we
state some still-open problems in the area and present some computer-obtained results.

\section{Preliminaries}

We begin by introducing some conventions and definitions.  In this, we mostly follow the earlier
work of Ball et~al. \cite{ball:rainbow}, where they introduced the following definitions, as well as
proved Lemma~\ref{lemma:obvious} (and its corollary) and Claim~\ref{claim:even}.

\begin{definition*}
  In the context of this paper, a \emph{coloring} is an edge-coloring of an undirected complete
  graph; the colors may come from an arbitrary set and there is no restriction that the coloring be
  \emph{proper}.  As above, a \emph{rainbow $n$-cycle} (again, sometimes called \emph{colorful} or
  \emph{totally multicolored} in other sources) within a coloring is a cycle consisting of $n$
  distinct vertices, all of whose edges are colored with distinct colors; in this case, we will also
  say that the coloring \emph{contains} a rainbow $n$-cycle (that is, the coloring restricted to the
  edges of the cycle is a bijection).  As notation, we write $(v_1, \dotsc, v_n)$ for the cycle that
  visits vertices $v_1, \dotsc, v_n$ in order (and then returns to $v_1$).  Notice that although we
  allow infinite graphs, all cycles are of course finite.
\end{definition*}

A simple lemma and an immediate corollary guide us in our study of rainbow cycles.

\begin{figure}[ht]
  \begin{pspicture}(-1.5,-1.5)(1.5,1.5)
    \SpecialCoor
    \pspolygon[showpoints=true]
    (1.5;15)  (1.5;45) (1.5;75) 
    (1.5;105)(1.5;135)(1.5;165)
    (1.5;195)(1.5;225)(1.5;255)
    (1.5;285)(1.5;315)(1.5;345)
    \psline[linestyle=dotted](1.5;15)(1.5;165)
    \rput( 0.9;90){$n$}
    \rput(-0.5;90){$m$}
    \psarc[arrows=<-]( 0.9;90){0.3}{270}{180}
    \psarc[arrows=<-](-0.5;90){0.3}{270}{180}
  \end{pspicture}
  \caption{An $(n+m-2)$-cycle divided into an $n$-cycle and an $m$-cycle by a
    chord. (See Lemma~\ref{lemma:obvious}.)}
  \label{figure:obvious}
\end{figure}
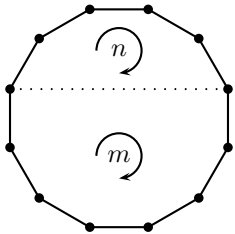

\begin{lemma} \label{lemma:obvious}
  If a coloring contains no rainbow $n$-cycles nor rainbow $m$-cycles, then it contains no rainbow
  $(n+m-2)$-cycles.
\end{lemma}
\begin{corollary} \label{corollary:obvious}
  In particular, if a coloring contains no rainbow $n$-cycle, then it contains no rainbow cycles of
  length $\ell$, where $\ell \equiv 2 \pmod{n-2}$.
\end{corollary}
\begin{proof}
  The proof is by contradiction: assume that the coloring contains a rainbow $(n+m-2)$-cycle.  It
  can be divided, by a single chord, into an $n$-cycle and an $m$-cycle (see
  Figure~\ref{figure:obvious}).  Consider the color of this chord.  On the one hand, it must agree
  with one of the other edges of the $n$-cycle to avoid a rainbow $n$-cycle; on the other hand, it
  must similarly agree with one of the other edges of the $m$-cycle to avoid a rainbow $m$-cycle.
  This, however, is a contradiction, as we assumed the outer $(n+m-2)$-cycle was rainbow.

  If an $n$-cycle is prohibited, then so are $(2n-2)$-cycles.  By induction, we obtain precisely all
  lengths congruent to $2\bmod n-2$.
\end{proof}

\section{Original Question and the Main Result}
\label{section:question}

The original question \cite{ball:rainbow} can be stated succinctly as follows:
\begin{question}
  Is the restriction in Lemma~\ref{lemma:obvious} the only restriction on what lengths of rainbow
  cycles a coloring can contain?  In other words, can any set of lengths that does not contradict
  the lemma be obtained as the lengths of rainbow cycles of some coloring?
\end{question}
\begin{remark*}
  As we shall see later in the proof of Lemma~\ref{lemma:monoid}, this question can be stated
  formally as follows: Is it true that for any submonoid $S\subset \N$ the set of lengths $\{ s+2
  \mid s\not\in S\}$, may be obtained as the lengths of rainbow cycles of some coloring?  Note also
  that if one only considers finite graphs, then it may be necessary to add in a restriction
  prohibiting sets containing arbitrarily large cycle lengths.
\end{remark*}

We may also formulate a weaker question based on Corollary~\ref{corollary:obvious}.

\begin{question} \label{question:weak}
  If $m \not\equiv 2 \pmod{n-2}$, does there exist a coloring of $K_m$ with a rainbow (Hamiltonian)
  $m$-cycle but no rainbow $n$-cycle?
\end{question}

However, both questions can be answered in the negative; indeed, they both contradict our main
theorem, which is somewhat of an opposite, Ramsey-type result.  We state it here but prove it later.

\newcounter{theorem:main}
\setcounter{theorem:main}{\value{theorem}}
\begin{theorem} \label{theorem:main}
  Suppose $n$ is an odd integer.  If a coloring does not contain a rainbow $n$-cycle, it also does
  not contain a rainbow $m$-cycle for all sufficiently large $m$; in particular, $m \ge 2n^2$
  suffices.
\end{theorem}

In particular, we shall show that if there is no rainbow $5$-cycle, there is also no rainbow
$10$-cycle, answering Question~\ref{question:weak} in the negative with $m=10$ and $n=5$.

\section{The Even Case}

Before proving our result, it is instructive to consider examples of colorings which contain some
lengths of rainbow cycles, but yet do not contain many other lengths.  In particular, we will
construct colorings that show that Theorem~\ref{theorem:main} is not true for even $n$.  Of the
following two results, Claim~\ref{claim:even} was proven in \cite{ball:rainbow} but
Claim~\ref{claim:2mod4} is original.

\begin{claim} \label{claim:even}
  There exists a coloring $\col$ of an infinite complete graph that contains rainbow cycles of all
  odd lengths, but no even-length rainbow cycles.  Furthermore, for each $n$, taking appropriate
  finite subgraphs (and their induced colorings) yields colorings that contain rainbow cycles of all
  odd lengths up to $n$, but still no even-length rainbow cycles.
\end{claim}
\begin{proof}
  We construct $\col$ first.  Let the vertex set be the positive integers $\Z^{\mathord +}$ and the
  colors the nonnegative integers $\N$ and define the color of the edge joining distinct vertices
  $x$ and $y$ to be
  \begin{displaymath}
    \col(x,y) = \begin{cases}
      0         & \text{if $y - x$ is even}, \\
      \min(x,y) & \text{if $y - x$ is odd}.
    \end{cases}
  \end{displaymath}
  First, we must show that there exist rainbow cycles of all odd lengths.  But this is easy!
  Consider the cycle $(1, 2, 3, \dotsc, k)$ for $k$ odd.  For $1 \le i < k$, the color of the edge
  joining $i$ and $i+1$ is $\col(i, i+1) = i$; finally, the color of the edge joining $k$ and $1$ is
  $\col(k,1) = 0$.  These are all distinct, so it remains to show that there are no even-length
  cycles.

  However, by Corollary~\ref{corollary:obvious}, we need only show that there are no rainbow
  $4$-cycles.  Suppose, by contradiction, that there is a rainbow cycle $(a, b, c, d)$.  How many
  times is the case ``$\col(x,y) = 0$ if $y - x$ is even'' used along the edges?  It cannot be used
  more than once because otherwise the cycle would contain a repeated color, $0$.  But $(b-a) +
  (c-b) + (d-c) + (a-d) = 0$ and thus, by parity, an even number of $b-a$, $c-b$, $d-c$, and $a-d$
  are odd and this case cannot be used exactly once.  Therefore, all of the edges use the
  ``$\col(x,y) = \min(x,y)$ if $y - x$ is odd'' case of the above definition.  Now, without loss of
  generality, assume $a$ is the smallest-numbered vertex of the four; a contradiction is immediate:
  since $a$ is the smallest-numbered vertex, $\col(a,b) = \col(a,d) = a$.  Therefore there are no
  rainbow $4$-cycles and thus no even-length rainbow cycles at all.

  Finally, taking the induced subgraph on the vertices from $1$ to $n$ accomplishes the second
  statement of the claim.  Indeed, all of the necessary odd-length rainbow cycles mentioned above
  still exist, and a rainbow $4$-cycle still doesn't exist.
\end{proof}

\begin{claim} \label{claim:2mod4}
  There exists another similar coloring $\col'$ of an infinite complete graph that contains rainbow
  cycles of all lengths, except those lengths congruent to $2\bmod 4$.  Furthermore, taking
  appropriate finite subgraphs has the same effect as before.
\end{claim}

\begin{proof}
  Again, let the vertex set be the positive integers $\Z^{\mathord +}$ and the colors the
  nonnegative integers $\N$; then for distinct vertices $x$ and $y$, define
  \begin{displaymath}
    \col'(x,y) = \begin{cases}
      0 & \text{if $y - x \equiv \hphantom{-}0 \pmod 2$}, \\
      x & \text{if $y - x \equiv \hphantom{-}1 \pmod 4$}, \\
      y & \text{if $y - x \equiv           - 1 \pmod 4$}.
    \end{cases}
  \end{displaymath}
  Note that $\col'$ is well-defined, as $\col'(x,y) = \col'(y,x)$ for all $x$ and $y$.  We show the
  existence of rainbow $k$-cycles, for all $k\not\equiv 2 \pmod 4$; luckily, the cycle $(1, 2, 3,
  \dotsc, k)$ accomplishes the task.  For $1 \le i < k$, the color of the edge joining $i$ and $i+1$
  is $\col'( i, i+1 ) = i$; we must only analyze the color of the edge joining $k$ and $1$.  If $k$
  is odd, then $\col'( k, 1 ) = 0$; if $k$ is divisible by $4$, the color $\col'( k, 1 ) = k$.  In
  either case, the cycle is rainbow.

  It remains to show that there is no $6$-cycle (once again, Corollary~\ref{corollary:obvious}
  implies we need only check this length).  We proceed by contradiction: assume there is a rainbow
  cycle $(a, b, c, d, e, f)$.  As before, the rule ``$\col'(x, y) = 0$ if $y - x \equiv 0 \pmod 2$''
  cannot be used at all.  We may also assume $b - a \equiv -1 \pmod 4$: if it is $1 \bmod 4$, simply
  change the direction of the cycle $(a, b, c, d, e, f) \mapsto (b, a, f, e, d, c)$.  It follows
  that $c - b \equiv -1$; indeed, if $c - b \equiv 1 \pmod 4$, then $\col( a , b ) = b = \col( b, c
  )$.  Similarly, we may conclude that $b - a \equiv c - b \equiv d - c \equiv \dotsb \equiv a - f
  \equiv -1 \pmod 4$.  This is a contradiction; clearly $(b - a) + (c - b) + \dotsb + (a - f) = 0$,
  but the former would imply it were equal to $2 \bmod 4$.

  Of course, taking the corresponding induced subgraphs achieves the finite results.
\end{proof}

\section{Building up to the Main Result}
\label{section:weaker}

We now prove a theorem slightly weaker than the main result.

\begin{theorem}[weaker version of Theorem~\ref{theorem:main}] \label{theorem:weaker}
  Suppose $n$ is an odd integer.  If a coloring does not contain a rainbow $n$-cycle, it also does
  not contain a rainbow $m$-cycle for all sufficiently large $m$; in particular, a \emph{cubic}
  bound $m \ge n^3/2$ suffices.
\end{theorem}

The proof will consist of two intermediate lemmas.

\begin{lemma} \label{lemma:hard}
  If a coloring does not contain a rainbow $n$-cycle, where $n = 2k+1$ is odd, it also does not
  contain a rainbow $m$-cycle, where $m = \binom{n}{2} = k \cdot (2k+1)$.
\end{lemma}
\begin{remark*}
  The $k=2$ case yields the result involving $5$ and $10$ mentioned in
  Section~\ref{section:question}.
\end{remark*}
\begin{proof}
  We prove the contrapositive.  Assume that we have a coloring, $\col$, of $K_m$ such that there is
  a rainbow (Hamiltonian) $m$-cycle but no rainbow $n$-cycle.  Without loss of generality, we may
  number the vertices of the $K_m$ by residues modulo~$m$ and insist that $\col( i, i+1 ) = i$ for
  $i \in \Z/m\Z$; in other words, we assume that our rainbow $m$-cycle is given by $(0, 1, \dotsc,
  m-1)$.  Consider the following $k+1$ different $(2k+1)$-cycles: (see Figure~\ref{figure:hard} for
  a visual accompaniment)
  \begin{displaymath}
    \begin{array}{@{($ $}r@{,$ $}r@{,$ $}r@{,$ $\dotsc,$ $}r@{}l@{$ $)}}
      0 & 1 & 2 & 2k &\\
      2k & 2k+1 & 2k+2 & 2\cdot 2k & \\
      \multicolumn{5}{c}{\vdots} \\
      (k-1)\cdot 2k & (k-1)\cdot 2k+1 & (k-1)\cdot 2k+2 & k\cdot
      2k & \\
      k\cdot 2k & k\cdot 2k+1 & k\cdot 2k+2, & (k+1)\cdot 2k
      &{}\equiv k
    \end{array}
  \end{displaymath}

  \begin{figure}[ht]
    \begin{pspicture}(-1.5,-1.5)(1.5,1.5)
      \SpecialCoor
      \pscircle(0;0){1.5}
      \psline[showpoints=true](1.5;110)(1.5;30)(1.5;-50)
      \psline[showpoints=true](1.5;70)(1.5;150)(1.5;230)
      \psline[linestyle=dotted](1.5;-50)(1.5;230)
      \uput[110](1.5;110){$0$}
      \uput[70](1.5;70){$k$}
      \uput[30](1.5;30){$2k$}
      \uput[-50](1.5;-50){$4k$}
      \uput[150](1.5;150){$-k$}
      \uput[230](1.5;230){$-3k$}
    \end{pspicture}
    \caption{The relevant chords of Lemma~\ref{lemma:hard}.}
    \label{figure:hard}
  \end{figure}
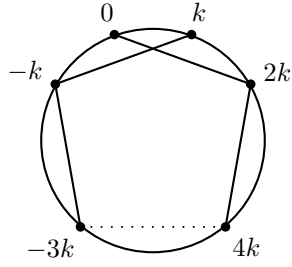

  By assumption, each of them must have a repeated color, so defining $c_i = \col( i\cdot k,
  (i+2)\cdot k)$ for $i\in \Z/m\Z$, we may conclude that
  \begin{displaymath}
    \begin{array}{r@{=\col($ $}r@{,$ $}r@{$ $)\in\{$ $}r@{,$ $}r@{,$ $\dotsc,$ $}r@{\}}l}
      c_0 & 0 & 2k & 0 & 1 & 2k-1 & , \\
      c_2 & 2k & 2\cdot 2k & 2k & 2k+1 & 2\cdot 2k - 1 & , \\
      \multicolumn{7}{c}{\vdots} \\
      c_{2(k-1)} & (k-1) \cdot 2k & k\cdot 2k & (k-1)\cdot 2k & (k-1)\cdot
      2k+1 & k\cdot 2k - 1 & , \\
      c_{-1} = c_{2k} & k \cdot 2k & k & k\cdot 2k & k\cdot 2k+1 & k - 1 & .
    \end{array}
  \end{displaymath}
  Now consider the $(2k+1)$-cycle $(0, 2k, 4k, \dotsc, (k-1)\cdot 2k, k\cdot 2k, k, k-1, k-2,
  \dotsc, 1)$; it has colors $\{0, 1, \dotsc, k-1\}\cup \{ c_0, c_2, \dotsc, c_{2k}\}$.  This
  collection must have a repeated color.  But none of $c_2, \dotsc, c_{2(k-1)}$ can contribute a
  repeated color, so we can conclude that one of $c_0$ and $c_{-1} = c_{2k}$ is a member of $\{0, 1,
  \dotsc, k-1\}$.  Notice, now, that we may translate this argument to also conclude that ``either
  $c_i$ or $c_{i-1}$ is a member of $\{i \cdot k, i\cdot k+1, \dotsc, (i+1)\cdot k - 1\}$''.

  By symmetry, assume that $c_0 \in \{0, 1, \dotsc, k-1\}$.  It follows that $c_1 \in \{k, k+1,
  \dotsc, 2k-1\}$ and, in general, $c_i \in \{i\cdot k, i\cdot k+1, \dotsc, (i+1)\cdot k - 1\}$.
  Now for the contradiction: consider the cycle $(0, 2k, \dotsc, k\cdot 2k, k, 3k, \dotsc,
  k\cdot(2k-1))$, illustrated in Figure~\ref{figure:star}.  The colors along the edges are precisely
  the $c_i$, no two of which can be equal! \qedhere

  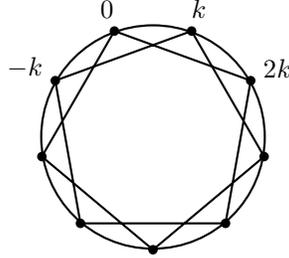
\begin{figure}[ht]
    \begin{pspicture}(-1.5,-1.5)(1.5,1.5)
      \SpecialCoor
      \pscircle(0;0){1.5}
      \pspolygon[showpoints=true]
      (1.5;30) (1.5;110)(1.5;190)
      (1.5;270)(1.5;350)(1.5;70) 
      (1.5;150)(1.5;230)(1.5;310)
      \uput[110](1.5;110){$0$}
      \uput[70](1.5;70){$k$}
      \uput[30](1.5;30){$2k$}
      \uput[150](1.5;150){$-k$}
    \end{pspicture}
    \caption{The eventual contradiction in the proof of Lemma~\ref{lemma:hard}.}
    \label{figure:star}
  \end{figure}
\end{proof}

\begin{lemma} \label{lemma:monoid}
  If a coloring contains no rainbow $n$-cycles, where $n = 2k+1$, and no rainbow $m$-cycles, where
  $m = k\cdot(2k+1)$, then it contains no rainbow $M$-cycles for all $M\ge 4k^3 - 2k^2 - 8k + 8 =
  n^3/2 - 2n^2 - 3n/2 + 11$.
\end{lemma}
\begin{proof}
  Let $\N(2)$ denote the set $\{2, 3, 4, \dotsc\}$.  With the operation $m\circ n = m + n - 2$ of
  Lemma~\ref{lemma:obvious}, $\N(2)$ becomes a monoid; moreover, by the map $n \mapsto n-2$, it is
  isomorphic to $\N$ under addition.  This observation is useful because of its interaction with
  ``absent'' rainbow cycle lengths.
  \begin{definition*}
    Let the \emph{spectrum} of a coloring $\col$ be the set of absent lengths of rainbow cycles,
    that is, $\{ n\ge 2 \mid \col $ does not contain a rainbow $n$-cycle$ \}$.  Notice that by
    definition, all spectra contain $2$.
  \end{definition*}
  Then by Lemma~\ref{lemma:obvious}, the spectrum of a coloring is a submonoid of $\N(2)$.  This
  allows us to apply the well-known Claim~\ref{claim:monoid}, perhaps first published by Sylvester,
  to complete the proof.
  \begin{claim}[``Frobenius coin-exchange problem'' \cite{sylvester}] \label{claim:monoid}
    If a submonoid of $\N$ under addition contains the relatively prime integers $a$ and $b$, then
    it contains all integers greater than or equal to $(a-1) (b-1)$.  (In other words, every
    sufficiently large integer can be written as a nonnegative integer linear combination of $a$ and
    $b$.)
  \end{claim}
  By the isomorphism between $(\N(2),\mathord{\circ})$ and $(\N, \mathord{+})$, we can check that
  $(2k+1) - 2$ and $k\cdot(2k+1) - 2$ are relatively prime (indeed, this is true because
  $k\cdot(2k+1) - 2 \equiv -1 \pmod{(2k+1)-2}$, a fact we will use later as well) and obtain the
  desired bound $4k^3 - 2k^2 - 8k + 8$ after a basic computation.
\end{proof}

Clearly, Lemma~\ref{lemma:hard} and Lemma~\ref{lemma:monoid} together complete the proof of
Theorem~\ref{theorem:weaker}, since $n^3/2 \ge n^3/2 - 2n^2 - 3n/2 + 11$ for $n \ge 2$.

\section{Proof of the Main Result} \label{section:main}

We now prove the main result, which is restated here for convenience:

\newcounter{temporary}
\setcounter{temporary}{\value{theorem}}
\setcounter{theorem}{\value{theorem:main}}
\begin{theorem}
  Suppose $n$ is an odd integer.  If a coloring does not contain a rainbow $n$-cycle, it also does
  not contain a rainbow $m$-cycle for all sufficiently large $m$; in particular, $m \ge 2n^2$
  suffices.
\end{theorem}
\setcounter{theorem}{\value{temporary}}

As in Section~\ref{section:weaker}, the proof will consist of two intermediate lemmas.

\begin{lemma} \label{lemma:harder}
  If a coloring does not contain a rainbow $n$-cycle, where $n = 2k+1 > 3$ is odd, it also does not
  contain a rainbow $m$-cycle, where $m = 3n-6 = 6k-3$.
\end{lemma}
\begin{remark*}
  The $k=2$ case yields the result that the absence of rainbow $5$-cycles implies the absence of
  rainbow $9$-cycles.
\end{remark*}
\begin{proof}
  The proof of this lemma has the same basic structure as that of Lemma~\ref{lemma:hard}, but the
  cycles we consider are different.  As before, we prove the contrapositive.  Assume that we have a
  coloring $\col$ of $K_m$ such that there is a rainbow (Hamiltonian) $m$-cycle but no rainbow
  $n$-cycle.  Without loss of generality, we may number the vertices of the $K_m$ by residues
  modulo~$m$ and insist that $\col( i, i+1 ) = i$ for $i \in \Z/m\Z$.

  Consider the following three $n$-cycles: (see Figure~\ref{figure:harder1} for a visual
  accompaniment --- all of the chords of the cycles are either along drawn edges or along the arcs
  of the outside cycle)
  \begin{displaymath}
    \begin{array}{lrrrrrcrr}
      ( & 0, & 1, & 2, & 3, & 4, & \dotsc, & n-1 & )\\
      ( & 1, & 0, & -1, & -2, & -3, & \dotsc, & 2n-4 & ) \\
      ( & 1, & 0, & n-1, & n, & n+1, & \dotsc, & 2n-4 & )
    \end{array}
  \end{displaymath}
  
  \begin{figure}[ht]
    \begin{pspicture}(-1.5,-1.5)(1.5,1.5)
      \SpecialCoor
      \pscircle(0;0){1.5}
      \psline[showpoints=true](1.5;76.66)(1.5;223.33)
      \psline[showpoints=true](1.5;103.33)(1.5;316.66)
      \psdots(1.5;343.33)(1.5;196.66)
      \uput[103.33](1.5;103.33){$0$}
      \uput[76.66](1.5;76.66){$1$}
      \uput[343.33](1.5;343.33){$n-2$}
      \uput[316.66](1.5;316.66){$n-1$}
      \uput[223.33](1.5;223.33){$2n-4$}
      \uput[196.66](1.5;196.66){$2n-3$}
      \uput[30](0.430205;30){$x$}
      \uput[150](0.430205;150){$z$}
    \end{pspicture}
    \caption{The chords in the first step of Lemma~\ref{lemma:harder}, as well as useful vertices to
      consider immediately after.}
    \label{figure:harder1}
  \end{figure}
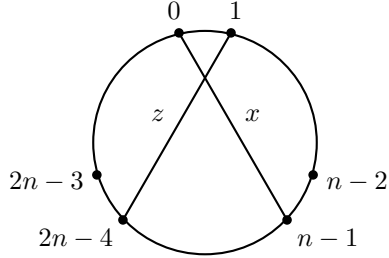

  Letting $x = \col(0, n-1)$ and $z = \col(1, 2n-4)$, the colors along the edges of these cycles are
  as follows:
  \begin{displaymath}
    \begin{array}{lrrrrrcrrr}
      ( & 0, & 1, & 2, & 3, & 4, & \dotsc, & n-2, & x )\\
      ( & 0, & -1, & -2, & -3, & -4, & \dotsc, & 2n-4, & z ) \\
      ( & 0, & x, & n-1, & n, & n+1, & \dotsc, & 2n-3, & z )
    \end{array}
  \end{displaymath}
  Because each of these cycles must contain a repeated color, it follows that either $x$ or $z$ is
  $0$.  (Assume there are no rainbow cycles.  Then if $x$ is not $0$, it is between $1$ and $n-2$;
  similarly, if $z$ is not $0$, it is between $-1$ and $-(n-2)$.  It follows that if neither $x$ nor
  $z$ is $0$, the third cycle is a rainbow cycle, a contradiction.)  Furthermore, let $y =
  \col(n-2,2n-3)$; by symmetry, we can also apply this argument to $x$ and $y$ to deduce that either
  $x$ or $y$ is $n-2$, and similarly for $y$ and $z$.  Altogether, either $x=0$, $y=n-2$, and
  $z=2n-4$ (call this ``orientation~$+$''); or $x=n-2$, $y=2n-4$, and $z=0$ (call this
  ``orientation~$-$'').

  In a sense, our previous argument was ``centered'' on the edge $(0,1)$.  By additive symmetry, we
  may also use the same argument centered at $(i,i+1)$ for any residue $i\in \Z /m\Z$; thus, we can
  assign an ``orientation~$\pm$'' to every residue $i$ (by construction, though, the orientations of
  $i$, $i + m/3$, and $i + 2m/3$ are the same).  Let $\Delta$ be a fixed integer; since $m$ is odd,
  there exists a residue $x$ such that $x$ and $x+\Delta$ have the same orientation.  Without loss
  of generality, we may assume that $x = 0$ and that the orientation is ``orientation~$+$.''  In the
  following, we let $\frac{n-3}2$ be our particular choice of $\Delta$ (note that since $n>3$,
  $\Delta > 0$).  Figure~\ref{figure:harder2} illustrates all of the edge colors that we may assume
  without any loss of generality.

  \renewcommand{\thesubfigure}{}
  \begin{figure}[ht]
    \subfigure[``Orientation~$+$'' centered at $(0,1)$]{
    \begin{pspicture}(-3,-2.25)(3,1.5)
      \SpecialCoor
      \pscircle(0;0){1.5}
      \psline[showpoints=true](1.5;76.66)(1.5;223.33)
      \psline[showpoints=true](1.5;103.33)(1.5;316.66)
      \psline[showpoints=true](1.5;343.33)(1.5;196.66)
      \uput[103.33](1.5;103.33){$0$}
      \uput[76.66](1.5;76.66){$1$}
      \uput[343.33](1.5;343.33){$n-2$}
      \uput[316.66](1.5;316.66){$n-1$}
      \uput[223.33](1.5;223.33){$2n-4$}
      \uput[196.66](1.5;196.66){$2n-3$}
      \uput[30](0.4;30){$0$}
      \uput[270](0.4;270){$n-2$}
      \uput[150](0.3;150){$2n-4$}
    \end{pspicture}
    }
    \subfigure[``Orientation~$+$'' centered at $(\Delta,\Delta+1)$ where $\Delta = \frac{n-3}2$]{
    \begin{pspicture}(-3,-2.25)(3,1.5)
      \SpecialCoor
      \pscircle(0;0){1.5}
      \psline[showpoints=true](1.5;283.33)(1.5;136.66)
      \psline[showpoints=true](1.5;256.66)(1.5;43.33)
      \psline[showpoints=true](1.5;16.66)(1.5;163.33)
      \rput[r](1.75;256.66){$\Delta+n-1$}
      \rput[l](1.75;283.33){$\Delta+n-2$}
      \uput[16.66](1.5;16.66){$\Delta+1$}
      \uput[43.33](1.5;43.33){$\Delta$}
      \uput[136.66](1.5;136.66){$\Delta+2n-3$}
      \uput[163.33](1.5;163.33){$\Delta+2n-4$}
      \uput[330](0.4;330){$\Delta$}
      \uput[90](0.4;90){$\Delta+2n-4$}
      \uput[210](0.3;210){$\substack{\textstyle\Delta+n \\ \textstyle -2}$}
    \end{pspicture}
    }
    \caption{The edge colors we may assume without loss of generality in the next step of
      Lemma~\ref{lemma:harder}.  They are illustrated on two separate diagrams for clarity.}
    \label{figure:harder2}
  \end{figure}
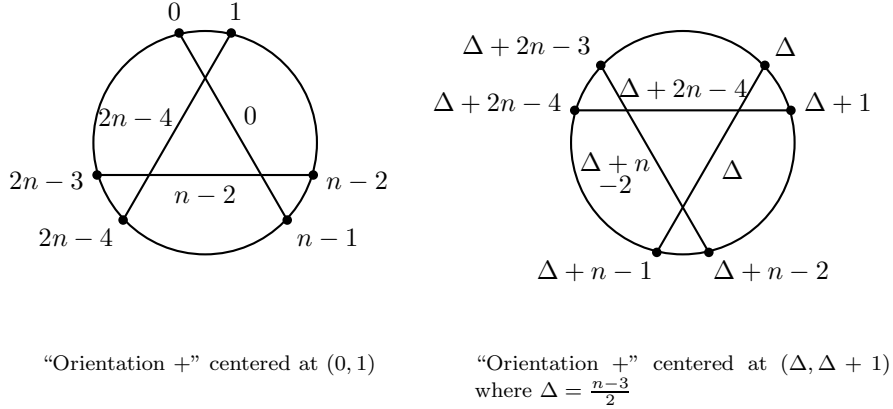

  Now let $t = \col(0,n-2)$; we will show that $t\in\{0, n-2\}$.  Before continuing, let us define
  the shorthand notation $x\upto y$ to mean $x, x+1, x+2, \dotsc, y$ (read: ``$x$ up to $y$'') and
  $x\dnto y$ to mean $x, x-1, x-2, \dotsc, y$ (read: ``$x$ down to $y$''); Then in this notation,
  examine the following $n$-cycles and their corresponding edge colors, illustrated separately in
  Figure~\ref{figure:harder3}:
  \begin{displaymath}
    \begin{array}{rll}
      & \multicolumn{2}{c}{\text{Respectively, the cycle and its edge colors}} \\
      \text{(a)} & (0\upto\Delta, \Delta+n-1\dnto n-2) & (0\upto \Delta, \Delta+n-2\dnto n-2, t) \\
      \text{(b)} & (0\upto\Delta+1, \Delta+2n-4\dnto 2n-3, n-2) &
      (0\upto \Delta, \Delta+2n-4\dnto 2n-3, n-2, t) \\
      \multirow{2}{*}{\text{(c)}} & (0,n-1\upto\Delta+n-2, & (0, n-1\upto \Delta+n-2, \\
      & \multicolumn{1}{r}{\Delta+2n-3\dnto 2n-3,n-2)}
      & \multicolumn{1}{r}{\Delta+2n-2\dnto 2n-3, n-2, t)}
    \end{array}
  \end{displaymath}

  \renewcommand{\thesubfigure}{(\alph{subfigure})}
  \begin{figure}[ht]
    \subfigure[$(0\upto\Delta, \Delta+n-1\dnto n-2)$]{
    \begin{pspicture}(-2,-2.25)(2.5,1.5)
      \SpecialCoor
      \pscircle(0;0){1.5}
      \pspolygon[showpoints=true](1.5;90)(1.5;330)
      (1.5;316.667)(1.5;303.333)(1.5;290)(1.5;276.667)(1.5;263.333)
      (1.5;36.66)(1.5;50)(1.5;63.33)(1.5;76.66)
      \uput[90](1.5;90){$0$}
      \uput[36.66](1.5;36.66){$\Delta$}
      \uput[263.33](1.5;263.33){$\Delta+n-1$}
      \uput[330](1.5;330){$n-2$}
      \uput[70](0.4;70){$t$}
      \uput[270](0.4;270){$\Delta$}
    \end{pspicture}
    }
    \subfigure[$(0\upto\Delta+1, \Delta+2n-4\dnto 2n-3, n-2)$]{
    \begin{pspicture}(-2.5,-2.25)(2.5,1.5)
      \SpecialCoor
      \pscircle(0;0){1.5}
      \pspolygon[showpoints=true]
      (1.5;90)(1.5;76.66)(1.5;63.33)(1.5;50)(1.5;36.66)(1.5;23.33)(1.5;156.66)
      (1.5;170)(1.5;183.33)(1.5;196.66)(1.5;330)
      \uput[90](1.5;90){$0$}
      \uput[20](1.4;20){$\Delta+1$}
      \uput[156.66](1.5;156.66){$\Delta+2n-4$}
      \uput[330](1.5;330){$n-2$}
      \uput[196.66](1.5;196.66){$2n-3$}
      \uput[10](0.4;10){$t$}
      \uput[90](0.5;90){$\Delta+2n-4$}
      \uput[270](0.6;270){$n-2$}
    \end{pspicture}
    }
    \subfigure[$(0,n-1\upto\Delta+n-2,\Delta+2n-3\dnto 2n-3,n-2)$.]{
    \begin{pspicture}(-2.5,-2.25)(2,1.5)
      \SpecialCoor
      \pscircle(0;0){1.5}
      \pspolygon[showpoints=true]
      (1.5;90)(1.5;330)(1.5;196.66)(1.5;183.33)(1.5;170)
      (1.5;156.66)(1.5;143.33)(1.5;276.66)(1.5;290)(1.5;303.33)
      (1.5;316.66)
      \uput[90](1.5;90){$0$}
      \uput[143.33](1.5;143.33){$\Delta+2n-3$}
      \uput[276.66](1.5;276.66){$\Delta+n-2$}
      \uput[330](1.5;330){$n-2$}
      \uput[316.66](1.5;316.66){$n-1$}
      \uput[196.66](1.5;196.66){$2n-3$}
      \uput[350](0.2;0){$0$}
      \uput[10](0.9;10){$t$}
      \uput[270](0.6;270){$n-2$}
      \rput(0.6;135){$\substack{\textstyle\Delta+n\\ \textstyle -2}$}
    \end{pspicture}
    }
    \caption{The three relevant $n$-cycles in the middle of the argument in
      Lemma~\ref{lemma:harder}.}
    \label{figure:harder3}
  \end{figure}

  Once more, because each of these cycles must contain a repeated color, $t$ can be only $0$ or
  $n-2$, as desired.  (Simply note that one of these cycles would be a rainbow cycle unless $t$
  repeats a color found on all of them: either $0$ or $n-2$.)  That is, $t =
  \col(0,n-2)\in\{0,n-2\}$; by symmetry, it follows that $\col(n-2,2n-4) \in \{n-2,2n-4\}$ and
  $\col(2n-4,0) \in \{2n-4, 0\}$.  By enumerating the cases, one can check that there necessarily
  exists an $i$ such that $\col\bigl(i(n-2), (i+1)(n-2)\bigr) + n-2 = \col\bigl((i+1)(n-2),
  (i+2)(n-2)\bigr)$.  (See Figure~\ref{figure:harder4}, particularly subfigure~(a), for continued
  visual accompaniment.)  Our final contradiction will be to show that this is impossible; by
  symmetry, we may assume that this happens when $i = 0$.

  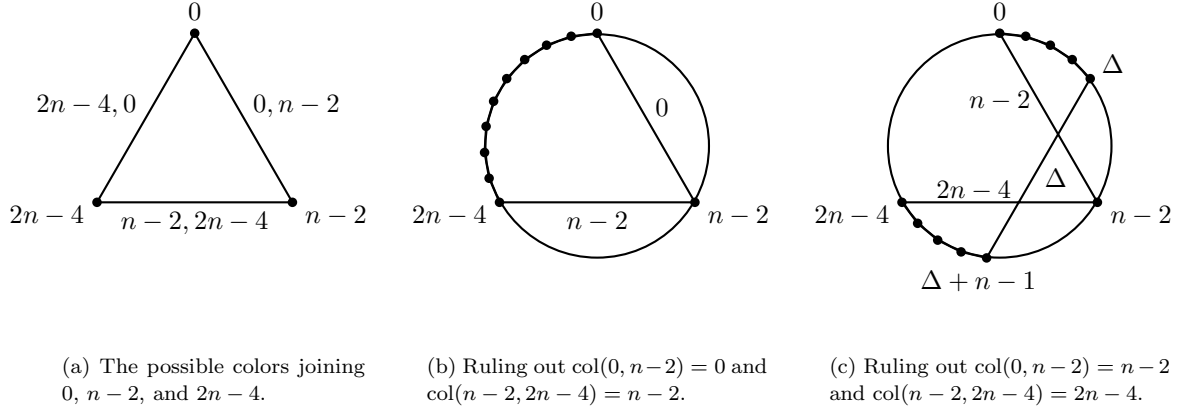
\begin{figure}[ht]
    \subfigure[The possible colors joining $0$, $n-2$, and $2n-4$.]{
    \begin{pspicture}(-2,-2.25)(2.5,1.5)
      \SpecialCoor
      \pspolygon[showpoints=true](1.5;90)(1.5;210)(1.5;330)
      \uput[90](1.5;90){$0$}
      \uput[330](1.5;330){$n-2$}
      \uput[210](1.5;210){$2n-4$}
      \uput[30](0.7;30){$0,n-2$}
      \uput[270](0.7;270){$n-2,2n-4$}
      \uput[150](0.7;150){$2n-4,0$}
    \end{pspicture}
    }
    \subfigure[Ruling out $\col(0,n-2)=0$ and $\col(n-2,2n-4)=n-2$.]{
    \begin{pspicture}(-2.5,-2.25)(2.5,1.5)
      \SpecialCoor
      \pscircle(0;0){1.5}
      \pspolygon[showpoints=true](1.5;90)(1.5;330)(1.5;210)
      (1.5;196.66)(1.5;183.33)(1.5;170)(1.5;156.66)(1.5;143.33)(1.5;130)(1.5;116.66)(1.5;103.33)
      \uput[90](1.5;90){$0$}
      \uput[330](1.5;330){$n-2$}
      \uput[210](1.5;210){$2n-4$}
      \uput[30](0.7;30){$0$}
      \uput[270](0.7;270){$n-2$}
    \end{pspicture}
    }
    \subfigure[Ruling out $\col(0,n-2)=n-2$ and $\col(n-2,2n-4)=2n-4$.]{
    \begin{pspicture}(-2.5,-2.25)(2.5,1.5)
      \SpecialCoor
      \pscircle(0;0){1.5}
      \pspolygon[showpoints=true](1.5;90)(1.5;330)(1.5;210)
      (1.5;223.33)(1.5;236.66)(1.5;250)(1.5;263.33)
      (1.5;36.66)(1.5;50)(1.5;63.33)(1.5;76.66)
      \uput[90](1.5;90){$0$}
      \uput[36.66](1.5;36.66){$\Delta$}
      \uput[263.33](1.5;263.33){$\Delta+n-1$}
      \uput[330](1.5;330){$n-2$}
      \uput[210](1.5;210){$2n-4$}
      \uput[90](0.3;90){$n-2$}
      \uput[250](0.3;250){$2n-4$}
      \uput[330](0.5;330){$\Delta$}
    \end{pspicture}
    }
    \caption{The eventual contradiction in Lemma~\ref{lemma:harder}.}
    \label{figure:harder4}
  \end{figure}

  We can rule out the possibility $\col(0,n-2) = 0$ and $\col(n-2,2n-4) = n-2$ by considering the
  $n$-cycle $(n-2,2n-4\upto 0)$ which would then have colors $(n-2, 2n-4\upto 0)$, a contradiction
  since no colors repeat and thus this cycle is rainbow.  We can rule out the only other possibility
  $\col(0,n-2) = n-2$ and $\col(n-2,2n-4) = 2n-4$ by considering the $n$-cycle $(0\upto\Delta,
  \Delta+n-1\upto 2n-4, n-2)$ which would then have colors $(0\upto \Delta, \Delta+n-1\upto 2n-4,
  n-2)$, also a contradiction.  We have thus obtained a contradiction in all cases!
\end{proof}

\begin{lemma} \label{lemma:adhoc}
  If a coloring contains no rainbow $n$-cycles, where $n = 2k+1 > 3$, and no rainbow $m_1$-cycles
  and $m_2$-cycles, where $m_1 = k\cdot(2k+1)$ and $m_2 = 6k-3$, then it contains no rainbow
  $M$-cycles for all $M\ge 8k^2 - 8k + 12 = 2n^2 - 13n + 23$.
\end{lemma}
\begin{proof}
  One approach is to simply proceed as in the proof of Lemma~\ref{lemma:monoid} and consider only
  $n$ and $m_2$ (since $n-2$ and $m_2-2$ are relatively prime).  This approach, using
  Claim~\ref{claim:monoid}, yields a bound of $12k^2 - 24k + 14 = 3n^2 - 18n + 29$, which is
  approximately a constant multiplicative factor worse than the desired bound.  We shall use both
  $m_1$ and $m_2$ to derive a better result.

  The analog of Claim~\ref{claim:monoid} for three variables is well-studied, and many special cases
  and algorithms have been developed.  For completeness, we use a self-contained ad-hoc argument
  applicable in our particular situation.

  \begin{claim} \label{claim:adhoc}
    Suppose a submonoid of $\N$ under addition contains the integers $a < b < c$, with $a \ge 3$
    odd, $b \equiv -2 \pmod a$, and $c \equiv -1 \pmod a$ (thus, $a$ and $b$ are relatively prime,
    as are $a$ and $c$).  Then the submonoid contains all integers greater than or equal to $ab/2 -
    a -3b/2 + c + 1$.
  \end{claim}
  \begin{proof}
    It suffices to prove the claim for the submonoid \emph{generated} by $a$, $b$, and $c$.

    For each residue $r$ modulo $a$, consider the smallest number in the submonoid congruent to $r$
    (modulo $a$).  The smallest number in this submonoid that is congruent to $r \equiv -(2x) \bmod
    a$, where $2x < a$, is $bx$; similarly, the smallest number congruent to $r \equiv -(2x+1) \bmod
    a$, where $2x+1 < a$, is $bx + c$.  The largest of these numbers, ranging over all residues $r$
    modulo $a$, occurs for $r\equiv 2 \equiv - (a-2) \equiv -(2 \frac{a-3}2 + 1) \bmod a$, in which
    case it is $y = b \frac{a-3}2 + c$.  By choice of $y$, all numbers greater than $y-a$ are in the
    monoid, as desired.
  \end{proof}

  To finish the proof of the lemma, apply this claim with $a = n-2$, $b = m_2-2$, and $c = m_1-2$,
  and use the isomorphism between $\N$ and $\N(2)$.  One obtains the desired bound of $8k^2 - 8k +
  12 = 2n^2 - 13n + 23$.
\end{proof}

As before, Lemma~\ref{lemma:harder} and Lemma~\ref{lemma:adhoc} together complete the proof of
Theorem~\ref{theorem:main}, since $2n^2 \ge 2n^2 - 13n + 23$ for $n \ge 2$.

\section{Further Directions}

This paper leaves open a few avenues of further experimentation which interest the author.  We state
some of these problems.
\begin{problem}
  Completely characterize when the existence of a rainbow $m$-cycle implies the existence of a
  rainbow $n$-cycle.
\end{problem}
\begin{problem}
  Let $g(n)$ be the smallest value of $M$ such that if a coloring does not contain a rainbow
  $n$-cycle, where $n$ is odd, then it also does not contain a rainbow $m$-cycle for all $m \ge M$.
  Determine $g(n)$ for specific cases or in general.
\end{problem}
\begin{remark*}
  For example, computer experimentation (see Appendix~\ref{appendix:computer}) yields $g(5) = 8$,
  $g(7) = 11$, $g(9) = 15$, and $15 \le g(11) \le 34$.
\end{remark*}
\begin{problem}
  In general, Theorem~\ref{theorem:main} shows that $g(n)$ is at most quadratic in $n$.  Is $g(n)$
  actually subquadratic?
\end{problem}
Finally, there is evidence to support the following conjecture.
\begin{conjecture*}
  The asymptotic behavior of a spectrum $S$ can be classified into three categories: either (a) $S$
  contains all sufficiently large numbers, (b) $S$ contains all sufficiently large even numbers, or
  (c) $S$ contains all sufficiently large numbers congruent to $2\bmod 4$.  In terms of monoids, the
  spectrum becomes regular either modulo one, two, or four.
\end{conjecture*}

\section*{Acknowledgments}

The author wishes to thank Professor Daniel Kleitman of the Massachusetts Institute of Technology
for guiding this work, as well as Jacob Fox, Matt Ince, and Petr Vojt\v{e}chovsk\'{y} for helpful
conversations.

\bibliographystyle{amsalpha}
%\bibliography{rainbow}
\newcommand{\etalchar}[1]{$^{#1}$}
\providecommand{\bysame}{\leavevmode\hbox to3em{\hrulefill}\thinspace}
\providecommand{\MR}{\relax\ifhmode\unskip\space\fi MR }
\providecommand{\MRhref}[2]{%
  \href{http://www.ams.org/mathscinet-getitem?mr=#1}{#2}
}
\providecommand{\href}[2]{#2}

\appendix
\section{Computer Results}
\label{appendix:computer}

An original C++ program was used to exhaustively determine (by searching through the space of all
possible colorings by always coloring the most-constrained edge), for small $n$ and $m$, whether or
not it is possible to color $K_m$ such that there exists a rainbow $m$-cycle but no rainbow
$n$-cycle.  These results are summarized in Table~\ref{table:computer}.  Some of these results have
been confirmed by Petr Vojt\v{e}chovsk\'{y}; in particular, the case $m = n+1 > 3$ is known (this
corresponds to the top entry of each column). \cite{ball:rainbow}

\newcommand{\easyno} {\textsf{x}}
\newcommand{\hardno} {\textsf{\textbf{X}}}
\newcommand{\easyyes}{\textsf{o}}
\newcommand{\hardyes}{\textsf{\textbf{O}}}
\newcommand{\toprow}[1]{\multicolumn{1}{c}{\textbf{#1}}}
\newcommand{\bottomrow}[1]{\multicolumn{1}{p{12pt}}{\centering\textbf{#1}}}
\begin{table}[p]
  \begin{tabular}{r|c|c|c|c|c|c|c|c|c|c|c|c|c|c|c|l}
     \multicolumn{1}{c}{}&\toprow{3} \\ \cline{2-2}
     \textbf{4} &\easyno&\toprow{4}  \\ \cline{2-3}
     \textbf{5} &\easyno&\easyyes&\toprow{5} \\ \cline{2-4}
     \textbf{6} &\easyno&\easyno&\hardyes&\toprow{6} \\ \cline{2-5}
     \textbf{7} &\easyno&\easyyes&\hardyes&\easyyes&\toprow{7} \\ \cline{2-6}
     \textbf{8} &\easyno&\easyno&\easyno&\easyyes&\hardyes&\toprow{8} \\ \cline{2-7}
     \textbf{9} &\easyno&\easyyes&\hardno&\easyyes&\hardyes&\easyyes&\toprow{9} \\ \cline{2-8}
     \textbf{10}&\easyno&\easyno&\hardno&\easyno&\hardyes&\hardyes&\hardyes&\toprow{10} \\ \cline{2-9}
     \textbf{11}&\easyno&\easyyes&\easyno&\easyyes&\hardno&\easyyes&\hardyes&\easyyes&\toprow{11} \\ \cline{2-10}
     \textbf{12}&\easyno&\easyno&\easyno&\easyyes&\easyno&\hardno&\hardyes&\easyyes&\hardyes&\toprow{12} \\ \cline{2-11}
     \textbf{13}&\easyno&\easyyes&\easyno&\easyyes&\hardno&\easyyes&\hardno&\easyyes&\hardyes&\easyyes&\toprow{13} \\ \cline{2-12}
     \textbf{14}&\easyno&\easyno&\easyno&\easyno&\hardno&\easyno&\hardyes&\hardyes&\hardyes& &\hardyes&\toprow{14} \\ \cline{2-13}
     \textbf{15}&\easyno&\easyyes&\easyno&\easyyes&\hardno&\easyyes&\hardno&\easyyes& &\easyyes& &\easyyes&\toprow{15} \\ \cline{2-14}
     \textbf{16}&\easyno&\easyno&\easyno&\easyyes&\easyno&\hardno&\easyno&\easyyes& & &\hardyes&\easyyes&\hardyes&\toprow{16} \\ \cline{2-15}
     \textbf{17}&\easyno&\easyyes&\easyno&\easyyes&\easyno&\easyyes&\hardno&\easyyes& &\easyyes& &\easyyes& &\easyyes&\toprow{17} \\ \cline{2-16}
     \textbf{18}&\easyno&\easyno&\easyno&\easyno&\easyno&\easyno&\hardno&\easyno& & & & & & &\hardyes&\textbf{18} \\ \cline{2-16}
     \textbf{19}&\easyno&\easyyes&\easyno&\easyyes&\easyno&\easyyes&\hardno&\easyyes&\hardno&\easyyes& &\easyyes& &\easyyes& &\textbf{19} \\ \cline{2-16}
     \textbf{20}&\easyno&\easyno&\easyno&\easyyes&\easyno&\easyno&\easyno&\easyyes&\easyno&\hardno& &\easyyes& & & &\textbf{20} \\ \cline{2-16}
     \textbf{21}&\easyno&\easyyes&\easyno&\easyyes&\easyno&\easyyes&\hardno&\easyyes& &\easyyes& &\easyyes& &\easyyes& &\textbf{21} \\ \cline{2-16}
     \textbf{22}&\easyno&\easyno&\easyno&\easyno&\easyno&\easyno&\easyno&\hardno&\hardno&\easyno& & & & & &\textbf{22} \\ \cline{2-16}
     \textbf{23}&\easyno&\easyyes&\easyno&\easyyes&\easyno&\easyyes&\easyno&\easyyes&\hardno&\easyyes& &\easyyes& &\easyyes& &\textbf{23} \\ \cline{2-16}
     \textbf{24}&\easyno&\easyno&\easyno&\easyyes&\easyno&\easyno&\easyno&\easyyes& &\hardno&\easyno&\easyyes& & & &\textbf{24} \\ \cline{2-16}
     \textbf{25}&\easyno&\easyyes&\easyno&\easyyes&\easyno&\easyyes&\easyno&\easyyes&\hardno&\easyyes& &\easyyes& &\easyyes& &\textbf{25} \\ \cline{2-16}
     \textbf{26}&\easyno&\easyno&\easyno&\easyno&\easyno&\easyno&\easyno&\easyno&\hardno&\hardno& &\easyno& & & &\textbf{26} \\ \cline{2-16}
     \textbf{27}&\easyno&\easyyes&\easyno&\easyyes&\easyno&\easyyes&\easyno&\easyyes& &\easyyes& &\easyyes& &\easyyes& &\textbf{27} \\ \cline{2-16}
     \textbf{28}&\easyno&\easyno&\easyno&\easyyes&\easyno&\easyno&\easyno&\easyyes&\easyno&\hardno& &\easyyes&\easyno& & &\textbf{28} \\ \cline{2-16}
     \textbf{29}&\easyno&\easyyes&\easyno&\easyyes&\easyno&\easyyes&\easyno&\easyyes&\easyno&\easyyes& &\easyyes& &\easyyes& &\textbf{29} \\ \cline{2-16}
     \textbf{30}&\easyno&\easyno&\easyno&\easyno&\easyno&\easyno&\easyno&\easyno& &\easyno& & & &\easyno& &\textbf{30} \\ \cline{2-16}
     \textbf{31}&\easyno&\easyyes&\easyno&\easyyes&\easyno&\easyyes&\easyno&\easyyes&\easyno&\easyyes& &\easyyes& &\easyyes& &\textbf{31} \\ \cline{2-16}
     \textbf{32}&\easyno&\easyno&\easyno&\easyyes&\easyno&\easyno&\easyno&\easyyes&\easyno&\easyno& &\easyyes& & &\easyno&\textbf{32} \\ \cline{2-16}
     \textbf{33}&\easyno&\easyyes&\easyno&\easyyes&\easyno&\easyyes&\easyno&\easyyes& &\easyyes& &\easyyes& &\easyyes& &\textbf{33} \\ \cline{2-16}
     \textbf{34}&\easyno&\easyno&\easyno&\easyno&\easyno&\easyno&\easyno&\easyno&\easyno&\easyno&\hardno&\hardno& & & &\textbf{34} \\ \cline{2-16}
     \textbf{35}&\easyno&\easyyes&\easyno&\easyyes&\easyno&\easyyes&\easyno&\easyyes&\easyno&\easyyes&\easyno&\easyyes& &\easyyes& &\textbf{35} \\ \cline{2-16}
     \textbf{36}&\easyno&\easyno&\easyno&\easyyes&\easyno&\easyno&\easyno&\easyyes&\easyno&\easyno& &\easyyes& & & &\textbf{36} \\ \cline{2-16}
     \textbf{37}&\easyno&\easyyes&\easyno&\easyyes&\easyno&\easyyes&\easyno&\easyyes&\easyno&\easyyes& &\easyyes& &\easyyes& &\textbf{37} \\ \cline{2-16}
     \textbf{38}&\easyno&\easyno&\easyno&\easyno&\easyno&\easyno&\easyno&\easyno&\easyno&\easyno& &\easyno& & & &\textbf{38} \\ \cline{2-16}
     \textbf{39}&\easyno&\easyyes&\easyno&\easyyes&\easyno&\easyyes&\easyno&\easyyes&\easyno&\easyyes& &\easyyes& &\easyyes& &\textbf{39} \\ \cline{2-16}
     \textbf{40}&\easyno&\easyno&\easyno&\easyyes&\easyno&\easyno&\easyno&\easyyes&\easyno&\easyno& &\easyyes&\hardno& & &\textbf{40} \\ \cline{2-16}
     \multicolumn{1}{c}{}
		 &\bottomrow{3}&\bottomrow{4}&\bottomrow{5}&\bottomrow{6}&\bottomrow{7}&\bottomrow{8}&\bottomrow{9}&\bottomrow{10}&\bottomrow{11}&\bottomrow{12}&\bottomrow{13}&\bottomrow{14}&\bottomrow{15}&\bottomrow{16}&\bottomrow{17}
  \end{tabular}
  \caption{Data obtained from a computer program used to determine when the absence of a particular
    length of rainbow cycles forbids the existence of another length.  For example, looking at
    column \textbf{5} and row \textbf{10} (the lowest `\hardno' in that column), we see that the
    absence of rainbow $5$-cycles implies the absence of rainbow $10$-cycles.  The meaning of the
    symbols is as follows: `\easyyes' means that Claim~\ref{claim:even} or Claim~\ref{claim:2mod4}
    implies that the existence of the larger length is possible, while `\hardyes' means that these
    Claims do not apply but the program produced a specific example.  `\hardno' means that the
    program determined that the larger length could not occur, while `\easyno's are consequences of
    Lemma~\ref{lemma:obvious} and earlier `\easyno's and `\hardno's.  An empty square means no
    results were obtained.}
  \label{table:computer}
\end{table}

\end{document}